\documentclass[11pt,a4paper]{article}
\usepackage[utf8]{inputenc}
\usepackage{geometry}
\usepackage[T1]{fontenc}
\usepackage{lmodern}
\usepackage{amsfonts}
\usepackage{amsmath}
\usepackage{amssymb}
\usepackage{amsthm}
\usepackage{mathtools}
\usepackage{fancyhdr}
\usepackage{enumitem}
\usepackage{mathrsfs}
\usepackage{hyperref}
\usepackage{cleveref}
\usepackage[all,cmtip]{xy}
\usepackage{manfnt}
\usepackage{marginnote}
\usepackage{tikz} 
\usepackage{tikz-cd}
\usepackage{graphicx}
\usepackage{mathabx}
\usepackage{tcolorbox}
\usepackage[mathscr]{euscript}
\usepackage{amsbsy}
\usepackage{float}

\usepackage[style = alphabetic]{biblatex}
\newtheoremstyle{mythm} 
    {.5em}                    
    {.5em}                    
    {\upshape}                   
    {}                           
    {\bfseries}                   
    {.}                          
    {.5em}                       
    {}  

\theoremstyle{mythm}
\newtheorem{thm}{Theorem}

\newtheorem*{notation}{Notation}



\usetikzlibrary{decorations.markings}

\makeatletter
\tikzcdset{
  open/.code     = {\tikzcdset{hook, circled};},
  closed/.code   = {\tikzcdset{hook, slashed};},
  open'/.code    = {\tikzcdset{hook', circled};},
  closed'/.code  = {\tikzcdset{hook', slashed};},
  circled/.code  = {\tikzcdset{markwith = {\draw (0,0) circle (.375ex);}};},
  slashed/.code  = {\tikzcdset{markwith = {\draw[-] (-.4ex,-.4ex) -- (.4ex,.4ex);}};},
  markwith/.code ={
    \pgfutil@ifundefined%
    {tikz@library@decorations.markings@loaded}%
    {\pgfutil@packageerror{tikz-cd}{You need to say %
      \string\usetikzlibrary{decorations.markings} to use arrows with markings}{}}{}%
    \pgfkeysalso{/tikz/postaction = {
      /tikz/decorate,
      /tikz/decoration={markings, mark = at position 0.5 with {#1}}}
    }
  },
}
\makeatother

\makeindex

\DeclareMathOperator{\Hom}{Hom}

\renewcommand{\S}{\text{S}}

\let \H \relax
\DeclareMathOperator{\H}{H}
\DeclareMathOperator{\C}{C}

\tikzset{commutative diagrams/.cd,
comm/.style={start anchor=center,end anchor=center,draw=none}
}

\title{Trivial Massey Product\\in Bounded Cohomology}
\author{Domenico Marasco}
\date{}

\begin{document}
\maketitle
\begin{abstract}
We show that in the bounded cohomology of non-abelian free groups the Massey triple product is always trivial when the second factor is represented by the coboundary of a decomposable quasi-morphism.

We also show that in the bounded cohomology of a negatively curved compact Riemannian manifold the Massey triple product is always trivial when the second factor is represented by an exact differential form.
\end{abstract}
\section{Introduction}
Bounded cohomology of groups has found various applications in different areas of mathematics, but computing it is a hard task, even in one of the most classic and studied cases: non-abelian free groups. Fix $n\geq2$ and denote by $F=F_n$ the non-abelian free group of rank $n$. Interestingly, the bounded cohomology of $F$ with trivial real coefficients $\H_b^k(F;\mathbb{R})$ is infinite dimensional when $k=2$ and $k=3$, while we still do not have any example of a non-trivial class for $k\geq4$. All classes in $\H_b^2(F;\mathbb{R})$ can be represented as the coboundaries of \textit{quasi-morphisms}, see Section 2.3.
There are various recent results investigating whether it is possible to construct a non-trivial bounded cocycle of degree $k\geq 4$ as the cup product of non-trivial quasi-morphisms; see \cite{BucherMonod+2018+1157+1162}, \cite{heuer2020cup}, \cite{fournierfacio2020infinite} and \cite{AmontovaBucher+2022+933+943}. All these results seem to suggest that $\cup\colon\H_b^2(F;\mathbb{R})\times\H_b^k(F;\mathbb{R})\to\H_b^{k+2}(F;\mathbb{R})$ could be trivial. In particular, in \cite{AmontovaBucher+2022+933+943} the authors generalize the results of the previous papers by showing that for every $\Delta$-decomposable quasi-morphism $\varphi$ and for every $\alpha\in\H_b^k(F;\mathbb{R})$, for $k\geq1$,
$$[\delta^1\varphi]\cup\alpha=0\in\H_b^{k+2}(F;\mathbb{R}).$$
A decomposition $\Delta$ provides a way to decompose each element of $F$ into finite subwords belonging to a fixed set of pieces. Furthermore, one requires $\Delta$ to behave well with respect to geodesic triangles in the Cayley graph, so that every bounded real function on the pieces induces a quasi-morphism on the free group. Quasi-morphisms constructed this way are called $\Delta$-decomposable. These concepts were first introduced by Heuer in \cite{heuer2020cup} and are defined precisely in Section 2.3.
It is important to note that the Brooks quasi-morphisms on non-selfoverlapping words and the Rolli quasi-morphisms are decomposable (Example 3.8 and Example 3.9 in \cite{heuer2020cup}).\\

In this paper we answer a question that naturally arises after the above results. In fact, when the cup product vanishes it is possible to define a three terms cohomology operation: the triple Massey product. More precisely, we can construct the product of three cohomology classes whenever the cup product of the first with the second and the cup product of the second with the third vanish (see Section 2.2 for the definition of Massey product). This means that in the context of the bounded cohomology of a free group it is possible to define the Massey product of three classes when the second term is represented by the coboundary of a $\Delta$-decomposable quasi-morphism. The following is the main result of the paper:
\begin{thm}
Let $\varphi$ be a $\Delta$-decomposable quasi-morphism, let $\alpha_1\in\H_b^{k_1}(F;\mathbb{R})$ and $\alpha_2\in\H_b^{k_2}(F;\mathbb{R})$, $k_1,k_2\geq1$. Then the Massey triple product $\langle\alpha_1,[\delta\varphi],\alpha_2\rangle$ in $\H_b^{k_1+k_2+1}(F;\mathbb{R})$ is trivial.
\end{thm}
This shows once more that finding a non-trivial class in $\H_b^k(F;\mathbb{R})$, for $k\geq4$, is challenging.\\

The Massey product has a connection with the LS-category. We list some results that point to a question connecting the Massey product in bounded cohomology and the ameanable category of a space.

The LS-category of a space $X$ is a homotopy invariant denoted by $cat_{LS}(X)$ and defined as the minimum number of (possibly disconnected) open subsets $\{U_i\}_{i\in I}$ necessary to cover $X$ so that each inclusion $U_i\hookrightarrow X$ is null-homotopic, for every $i\in I$. The \textit{cup lenght} of a space $cup(X)$ is the maximum number $k$ of singular cohomology classes $\alpha_1\dots,\alpha_k$ of positive degree such that $\alpha_1\cup\dots\cup\alpha_k\neq0$. It is a classic result that the cup length bounds from below the LS-category (\cite{cornea2003lusternik}, Proposition 1.5):
$$
cup(X)+1\leq cat_{LS}(X).
$$
Therefore, a space $X$ with $cat_{LS}(X)\leq2$ must have trivial cup product, which means that the Massey product is defined for any triple of elements in $\H^{\bullet}(X;\mathbb{R})$. Interestingly, we can say even more: Corollary 4.6 of \cite{RUDYAK199937} shows that
$$
cat_{LS}(X)\leq2\implies \langle\alpha_1,\alpha_2,\alpha_3\rangle\;\text{is trivial for every }\alpha_1,\alpha_2,\alpha_3\in\H^\bullet(X;\mathbb{R}).
$$
There is a generalization of the LS-category in the context of bounded cohomology. For a space $X$, its \textit{amenable category}, denoted by $cat_{Am}(X)$, is defined as the minimum number of open subsets $\{U_i\}_{i\in I}$ that cover $X$ such that the map induced by the inclusion $U_i\hookrightarrow X$ on the fundamental group $\pi_1(U_i,x)\to\pi_1(X,x)$ has amenable image, for every $x\in U_i$ and for every $i\in I$. In the context of bounded cohomology it is possible to define the \textit{bounded cup length} $cup_b(X)$ of a space $X$ exactly like its classic counter part. In Remark 3.17 of \cite{https://doi.org/10.48550/arxiv.2012.00612} the authors ask whether the following holds:
$$
cup_b(X)+1\leq cat_{Am}(X)?
$$
In particular, the fundamental group of a wedge of $n$ circles is $F$ and $cat_{Am}(\bigvee S^1)=2$, which means that a positive answer to the question above would show that $\cup\colon\H_b^\bullet(F;\mathbb{R})\times\H_b^\bullet(F;\mathbb{R})\to\H_b^{\bullet}(F;\mathbb{R})$ is always trivial in positive degrees. Given Theorem 1 and the connection between Massey product and LS-category it is natural to ask if the following could also hold:
$$
cat_{Am}(X)\leq2\implies \langle\alpha_1,\alpha_2,\alpha_3\rangle\;\text{is trivial for every }\alpha_1,\alpha_2,\alpha_3\in\H_b^\bullet(X;\mathbb{R})?
$$

We finally observe that Theorem 1 admits an analogous version in the context of bounded classes defined by differential forms.
Let $M$ be a compact Riemannian manifold with negative sectional curvatures and with (possibly empty) convex boundary. For every exact $2$-form $\xi=d\varphi\in\text{E}\Omega^2(M)$ there is a bounded cocycle $c_\xi$ defined by integrating $\xi$ over the straightening of simplices. These cocycles represent exact bounded cohomology classes $[c_\xi]\in\text{EH}_b^2(M)$. In \cite{https://doi.org/10.48550/arxiv.2202.04419} it is shown that the cup product of $[c_\xi]$ with any other bounded cohomology class $\alpha\in\H_b^k(M;\mathbb{R})$ is trivial, for $k\geq1$. This means that we can define the triple Massey product whenever the second term can be represented by a cocycle defined by an exact form. We will show the following:
\begin{thm}
Let $\xi\in\text{E}\Omega^2(M)$ be an exact $2$-form. Let $\alpha_1\in\H_b^{k_1}(M;\mathbb{R})$ and $\alpha_2\in\H_b^{k_2}(M;\mathbb{R})$, $k_1,k_2\geq1$. Then the Massey triple product $\langle\alpha_1,[c_\xi],\alpha_2\rangle$ in $\H_b^{k_1+k_2+1}(M;\mathbb{R})$ is trivial.
\end{thm}
\section*{Acknowledgements}
I thank Marco Moraschini who has suggested the topic of this paper. I also thank Pietro Capovilla for useful discussions about the amenable category.
\section{Preliminaries}
\subsection{Bounded Cohomology of Groups}
Let $G$ be a group, for every $k\in\mathbb{N}$ we set 
$
\C^k(G;\mathbb{R})=\{\varphi\colon G^k\to;\mathbb{R}\}.
$
The coboundary maps $\delta^k\colon\C^k(G;\mathbb{R})\to\C^{k+1}(G;\mathbb{R})$ are defined as follows:
\begin{align*}
\delta^{k}\varphi(g_1,\dots,g_{k+1})&=\varphi(g_2,\dots,g_{k+1})\\&+
\sum_{i=1}^k(-1)^i\varphi(g_1,\dots,g_i g_{i+1},\dots,g_{k+1})\\&+
(-1)^{k+1}\varphi(g_1,\dots,g_k).
\end{align*}
The cochain complex $(\C^\bullet(G;\mathbb{R}),\delta^\bullet)$ is the inhomogeneous resolution of $G$ with coefficients in $\mathbb{R}$ and its cohomology, denoted by $\H^\bullet(G;\mathbb{R})$, is the standard cohomology of $G$ with coefficients in $\mathbb{R}$. From now on we will omit the coefficients from the notation since we are only going to use trivial real coefficients.

We put the sup norm on $\C^k(G)$ and denote by $\C_b^k(G)$ the subspace of cochains with bounded norm. Since the coboundary of a bounded cochain is also bounded, $(\C^\bullet_b(G),\delta^\bullet)$ is a cochain complex and its homology $\H_b^\bullet(G)$ is called the \textit{bounded cohomology of }$G$. The map induced on cohomology by the inclusion $\C_b^\bullet(G)\subset\C^\bullet(G)$ is usually called the \textit{comparison map} $c_b^k\colon\H_b^k(G)\to\H^k(G).$
The kernel of $c_b^k$ is called the \textit{exact bounded cohomology of} $G$ and it is denoted by $\text{EH}_b^k(G)$.
\subsection{Massey Triple Product}
In this section we give the definition of the Massey triple product, a cohomology operation first defined by Massey in \cite{UeharaMassey+2015+361+377}. The classical context for this operation is singular cohomology, but the usual definition can be harmlessly transferred to the context of bounded cohomology.\\

Let $\alpha\in\H_b^{k}(G)$, $\alpha_1\in\H_b^{k_1}(G)$ and $\alpha_2\in\H_b^{k_2}(G)$. The Massey triple product of $\alpha_1$, $\alpha$, $\alpha_2$, denoted by $\langle\alpha_1,\alpha,\alpha_2\rangle$, can be defined whenever $\alpha_1\cup\alpha=0=\alpha\cup\alpha_2$. In general the output of this operation is not a single element, but a set of elements in $\H_b^{k_1+k+k_2-1}(G)$ constructed as follows. Pick a representative for each class $[\omega_1]=\alpha_1$, $[\omega]=\alpha$, $[\omega_2]=\alpha_2$
and choose $\beta_1\in\C_b^{k_1+k-1}(G)$ and $\beta_2\in\C_b^{k_2+k-1}(G)$ so that $\delta\beta_1=\omega_1\cup\omega$ and $\delta\beta_2=\omega\cup\omega_2$. The following class is in the Massey triple product:
$$
(-1)^k[(-1)^{k_1}\omega_1\cup\beta_2-\beta_1\cup\omega_2]\in\langle\alpha_1,\alpha,\alpha_2\rangle\subset\H_b^{k_1+k+k_2-1}(G).
$$
Observe that $(-1)^{k_1}\omega_1\cup\beta_2-\beta_1\cup\omega_2$ is a cocycle, in fact
$$
\delta((-1)^{k_1}\omega_1\cup\beta_2-\beta_1\cup\omega_2)
=\omega_1\cup\delta\beta_2-\delta\beta_1\cup\omega_2=0.
$$
If we choose different primitives for the cup products we could end up with a different element in the Massey product. In fact, let $\beta_i^\prime=\beta_i+\gamma_i$, with $\gamma_i\in\C_b^{k_i+k-1}(G)$ a cocycle (observe that all the possible primitives for $\omega\cup\omega_i$ are obtained this way). Computing the product with these primitives we get:
\begin{align*}
[(-1)^{k_1}\omega_1\cup\beta^\prime_2-\beta^\prime_1\cup\omega_2]&=
[(-1)^{k_1}\omega_1\cup(\beta_2+\gamma_2)-(\beta_1+\gamma_1)\cup\omega_2]\\&=
[(-1)^{k_1}\omega_1\cup\beta_2-\beta_1\cup\omega_2]
+[\omega_1]\cup(-1)^{k_1}[\gamma_2]-[\gamma_1]\cup[\omega_2].
\end{align*}
Since $\gamma_1$ and $\gamma_2$ are generic cocycles and the result only depends on their classes, we can more compactly describe the Massey triple product of $\alpha_1$, $\alpha$ and $\alpha_2$ as the set
$$
\langle\alpha_1,\alpha,\alpha_2\rangle=(-1)^k[(-1)^{k_1}\omega_1\cup\beta_2-\beta_1\cup\omega_2]
+\alpha_1\cup\H_b^{k+k_2-1}(G)+\H_b^{k_1+k-1}(G)\cup\alpha_2.
$$
It is easy to check that the set $\langle\alpha_1,\alpha,\alpha_2\rangle$ does not depend on the choice of the representatives for $\alpha_1$, $\alpha$ and $\alpha_2$.

This product is said to be trivial when $0\in\langle\alpha_1,\alpha,\alpha_2\rangle$. Observe that when the cup product is trivial the Massey product is a well defined operation, in fact in this case $\alpha_1\cup\H_b^{k+k_2-1}(G)=0=\H_b^{k_1+k-1}(G)\cup\alpha_2$ for every $\alpha_1$ and every $\alpha_2$.
\subsection{Quasi-morphisms and free groups}
A \textit{quasi-morphism} is a function $f\colon G\to\mathbb{R}$ that fails to be a homomorphism of groups by some bounded amount. More precisely, $f$ is a quasi-morphism if there is a constant $D\geq0$ such that for every $g_1,g_2\in G$,
$$
|f(g_1)|+|f(g_2)|-|f(g_1g_2)|<D.
$$
We denote by $\text{Q}(G)$ the set of quasi-morphisms of a group $G$. Clearly, homomorphisms of groups $\Hom(G;\mathbb{R})$ and bounded functions $\C_b^1(G)$ are quasi-morphisms, but functions that are the sum of a homomorphism and a bounded function are usually called \textit{trivial quasi-morphisms}. The study of the second exact bounded cohomology of $G$ is tied to the study of non-trivial quasi-morphisms. In fact, the coboundary of a quasi-morphism $f$ is a bounded $2$-cocycle $\delta f\in\C_b^2(G)$ and it comes directly from the definitions that the map $f\mapsto[\delta f]$ induces an isomorphism:
$$
\frac{\text{Q}(G)}{\C_b^1(G)\oplus\Hom(G;\mathbb{R})}\cong\text{EH}_b^2(G).
$$

Let $F=F_n$ be the non-abelian free group with $n$-generators, for $n\geq2$. Since $\H^2(F)=0$, the comparison map $c_b^2\colon\H_b^2(F)\to\H^2(F)$ is the zero map and $$\H_b^2(F)\cong\text{EH}_b^2(F)\cong\frac{\text{Q}(F)}{\C_b^1(F)\oplus\Hom(F;\mathbb{R})}.$$
This means that in order to compute the second bounded cohomology of $F$ it is enough to study its quasi-morphisms.
\subsection{Decompositions}
In this section we quickly overview decompositions of a non-abelian free group $F$ and decomposable quasi-morphisms. For more details we refer to Section 3 of \cite{heuer2020cup}, where Heuer first introduced these concepts.
\\

Let $\mathcal{P}\subset F$ be a symmetric subset with $1\notin\mathcal{P}$, we will call its elements \textit{pieces}. We denote by $\mathcal{P}^*$ the set of finite sequences of pieces, including the empty sequence. A decomposition of $F$ into the pieces $\mathcal{P}$ is a map $\Delta\colon F\to\mathcal{P}^*$ that sends every $g\in F$ to a finite sequence $\Delta(g)=(g^{(1)},\dots,g^{(m)})$ with $g^{(j)}\in\mathcal{P}$ so that the following properties are satisfied:
\begin{itemize}
    \item[(i)] $g=g^{(1)}\dots g^{(m)}$ as a reduced word. 
    \item[(ii)] $\Delta(g^{-1})=((g^{(m)})^{-1},\dots,(g^{(1)})^{-1})$.
    \item[(iii)] $\Delta(g^{(i)}\dots g^{(j)})=(g^{(i)},\dots,g^{(j)})$ for every $1\leq i\leq j \leq m$.
\item[(iv)]
Lastly, for every $g,h\in F$ the number of pieces of the thick part of the geodesic triangle $(1,g,gh)$ is uniformly bounded. More precisely, for every $g,h\in F$ there are $c_1$, $c_2$, $c_3$, $r_1$, $r_2$, $r_3\in F$ so that
\begin{itemize}
    \item[$\bullet$] $\Delta(g)=\Delta(c_1^{-1})\Delta(r_1)\Delta(c_2)$
    \item[$\bullet$] $\Delta(h)=\Delta(c_2^{-1})\Delta(r_2)\Delta(c_3)$
    \item[$\bullet$] $\Delta((gh)^{-1})=\Delta(c_3^{-1})\Delta(r_3)\Delta(c_1)$
\end{itemize}
where we take the $c_i$ so that $\Delta(c_i)$ is of maximal length. Now, there must be a constant $R>0$ such that the lenght of $\Delta(r_i)$ is smaller than $R$.
\end{itemize}
\begin{figure}[H]
    \centering
    \includegraphics[scale=0.8]{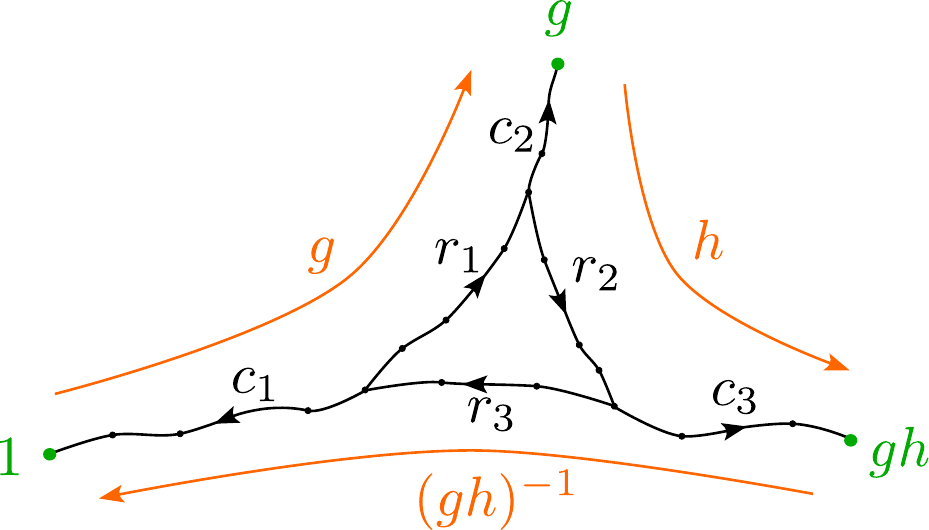}
    \caption{Decomposition of the geodesic triangle $(1,g,gh)$.}
\end{figure}

Let $\Delta$ be a decomposition with pieces $\mathcal{P}$. For $g\in G$ with $\Delta(g)=(g^{(1)},\dots,g^{(m)})$, we denote by $|\Delta(g)|=m$, the length of the decomposition of $g$. For any bounded alternating real function on the pieces $\lambda\in\ell^\infty_{alt}(\mathcal{P})$, we define $\varphi_{\lambda,\Delta}\colon F\to\mathbb{R}$ so that for every $g\in F$
$$
\varphi_{\lambda,\Delta}(g)=\sum_{j=1}^{|\Delta(g)|}\lambda(g^{(j)}).
$$
Any function constructed as above is called a $\Delta$\textit{-decomposable quasi-morphism}. Indeed, using property (iv) of decompositions it is straightforward to show that $\varphi_{\lambda,\Delta}$ is a quasi-morphism, see Proposition 3.6 of \cite{heuer2020cup}.

In Example 3.8 and Example 3.9 of \cite{heuer2020cup} the author shows that  the Brooks quasi-morphisms on non-selfoverlapping words and the Rolli quasi-morphisms are decomposable.
\subsection{Aligned Cochain Complex}
In this section we define the inhomogeneous version of the aligned cochain complex, which is another resolution that computes the bounded cohomology of $F$. The homogeneous version was first introduced in \cite{bucherbounded}. We will use the inhomogenous version of this resolution to simplify computations that involve decompositions, as it was done in \cite{AmontovaBucher+2022+933+943}.
\\

For every $k\geq0$ define 
$$
\mathcal{B}^k=\left\{(g_1,\dots,g_k)\in F^k\;\Big|\;g_ig_{i+1} \text{ is a reduced word, } g_i\neq 1\right\}.
$$
Define
$\mathcal{A}^k(F)=\{\varphi\colon\mathcal{B}^k\to\mathbb{R}\}$ and equip it with the usual inhomogeneous coboundary maps. We call  $(\mathcal{A}^\bullet(F),\delta^\bullet)$ the \textit{inhomogeneous aligned cochain complex}.
Denote by $\mathcal{A}_b^\bullet(F)$ the subcomplex of bounded cochains with respect to the sup norm. Furthermore, denote by $\mathcal{A}_{b,alt}^\bullet(F)$ the subcomplex of alternating cochains, i.e. of the cochains $\varphi\in\mathcal{A}_b^k(F)$ such that
$$
\varphi(g_1,\dots,g_k)=(-1)^{\lceil k/2\rceil}\varphi(g_k^{-1},\dots,g_1^{-1})
$$
for every $(g_1,\dots,g_k)\in\mathcal{B}^k$. There is a chain map $alt\colon\mathcal{A}_b^k(F)\to\mathcal{A}_{b,alt}^k(F)$ defined as
$$
alt(\varphi)(g_1,\dots,g_k)=\frac{1}{2}(\varphi(g_1,\dots,g_k)+(-1)^{\lceil k/2\rceil}\varphi(g_k^{-1},\dots,g_1^{-1}).
$$
The map above is the inhomogenous version of the usual alternating operator, which is well known to induce isomorphism on cohomology (see e.g. Theorem 1.4.2 in \cite{frigerio2018gromov}, Theorem 9.1 of \cite{eilenberg1944singular}, or Theorems 6.9 and 6.10 of \cite{eilenberg2015foundations}).
Furthermore, it is shown in Proposition 8 of \cite{BucherMonod+2018+1157+1162} that the composition of the restriction $r\colon\C_b^k(F)\to\mathcal{A}_b^k(F)$ with $alt$
$$
alt\circ r\colon\C_b^k(F)\to\mathcal{A}_{b,alt}^k(F)
$$
induces an isomophism between $\H_b^\bullet(F)$ and $\H^\bullet(\mathcal{A}_{b,alt}^k(F))$. This implies that the restriction $r\colon\C_b^k(F)\to\mathcal{A}_b^k(F)$ also induces an isomorphism on cohomology.
\section{Proof of Theorem 1}
Let $\varphi$ be a $\Delta$-decomposable quasi-morphism, let $\alpha_1\in\H_b^{k_1}(F)$, let $\alpha_2\in\H_b^{k_2}(F)$ and let $\omega_i\in\mathcal{A}_b^{k_i}(F)$ be the restriction of a representative for $\alpha_i$. In \cite{AmontovaBucher+2022+933+943} the authors show that $[\delta\varphi]\cup\alpha_i=0\in\H_b^{k_i+2}(F)$ by finding an explicit primitive $\beta_i\in\mathcal{A}_b^{k_i+1}(F)$ for $\delta\varphi\cup\omega_i\in\mathcal{A}_b^{k_i+2}(F)$. In what follows we will find a bounded primitive in $\mathcal{A}_b^{k_1+k_2}(F)$ for the cocycle
$$
(-1)^{k_1}\omega_1\cup\beta_2-\beta_1\cup\omega_2\in\mathcal{A}_b^{k_1+k_2+1}(F).
$$
This will be enough to show that $\langle\alpha_1,[\delta\varphi],\alpha_2\rangle$ is trivial because as we have seen in Section 2.5, the restriction $r\colon\C_b^\bullet(F)\to\mathcal{A}_b^\bullet(F)$ induces an isomorphism on the cohomology of the cochain complexes.

We write out the explicit formulas given in \cite{AmontovaBucher+2022+933+943} for the bounded primitives $\beta_1\in\mathcal{A}_b^{k_1+1}(F)$ and $\beta_2\in\mathcal{A}_b^{k_2+1}(F)$ of the cup products
$\omega_1\cup \delta\varphi=\delta\beta_1$ and $\delta\varphi\cup\omega_2=\delta\beta_2$:
$$
\beta_1=(-1)^{k_1}\omega_1\cup\varphi-\delta\eta_1
$$
$$
\beta_2=\varphi\cup\omega_2+\delta\eta_2,
$$
for suitably chosen $\eta_1\in\mathcal{A}^{k_1}(F)$ and $\eta_2\in\mathcal{A}^{k_2}(F)$.
In order to define $\eta_i\in\mathcal{A}^{k_i}(F)$ we first need to introduce the following notation. For any $g\in F$, if $\Delta(g)=(g^{(1)},\dots,g^{(N)})$ is its decomposition, then we denote by $z^<_j(g)$ and $z^>_j(g)$ the product of the first $(j-1)$, respectively the last $(N-j)$, pieces of the decomposition:
$$
z^<_j(g)=g^{(1)}g^{(2)}\dots g^{(j-1)}
$$
$$
z^>_j(g)=g^{(j+1)}g^{(j+2)}\dots g^{(N)}.
$$
Finally we define, for every $(g_1,\dots,g_{k_1})\in\mathcal{B}^{k_1}$,
\begin{align*}
\eta_1(g_1,\dots,g_{k_1})&=
\sum_{j=1}^{|\Delta(g_{k_1})|}
\omega_1(g_1,\dots,g_{k_1-1},z^<_j(g_{k_1}))\cdot\varphi(g_{k_1}^{(j)})
\end{align*}
and, for every $(h_1,\dots,h_{k_2})\in\mathcal{B}^{k_2}$,
\begin{align*}
\eta_2(h_1,\dots,h_{k_2})&=
\sum_{j=1}^{|\Delta(h_1)|}
\varphi(h_1^{(j)})\cdot\omega_2(z^>_j(h_1),h_2,\dots,h_{k_2}).
\end{align*}
Our goal is to find a bounded primitive for
\begin{align*}
(-1)^{k_1}\omega_1\cup\beta_2-\beta_1\cup\omega_2&=
(-1)^{k_1}\omega_1\cup(\varphi\cup\omega_2+\delta\eta_2)
-((-1)^{k_1}\omega_1\cup \varphi-\delta\eta_1)\cup\omega_2
\\&=(-1)^{k_1}\omega_1\cup\delta\eta_2+\delta\eta_1\cup\omega_2\in\mathcal{A}_b^{k_1+k_2+1}(F).
\end{align*}
Observe that the (not necessarily bounded) cochain $\omega_1\cup\eta_2+\eta_1\cup\omega_2\in\mathcal{A}^{k_1+k_2}(F)$ is a primitive. In order to find a bounded primitive we will define a cochain $\eta\in\mathcal{A}^{k_1+k_2-1}(F)$ such that $\omega_1\cup\eta_2+\eta_1\cup\omega_2-(-1)^{k_1}\delta\eta$ is bounded. Observe that this cochain is also a primitive, in fact
$$
\delta(\omega_1\cup\eta_2+\eta_1\cup\omega_2-(-1)^{k_1}\delta\eta)=\delta(\omega_1\cup\eta_2+\eta_1\cup\omega_2).
$$
Before defining $\eta$ we use the cocycle condition for $\omega_1$ and $\omega_2$ to expand $\omega_1\cup\eta_2$ and $\eta_1\cup\omega_2$, evaluated in $(g_1,\dots,g_{k_1},h_1,\dots,h_{k_2})\in \mathcal{B}^{k_1+k_2}$. In fact $\delta\omega_1=0$ implies that for every $1\leq j\leq|\Delta(h_1)|$,
\begin{align*}
\omega_1(g_1,\dots,g_{k_1})&=\omega_1(g_1,\dots,g_{k_1}z^<_j(h_1))\\&+
(-1)^{k_1}\sum_{i=1}^{k_1-1}(-1)^i\omega_1(g_1,\dots,g_i g_{i+1},\dots,g_{k_1},z^<_j(h_1))\\&-
(-1)^{k_1}\omega_1(g_2,\dots,g_{k_1},z^<_j(h_1)).
\end{align*}
Similarly, $\delta\omega_2=0$ implies that for every $1\leq j\leq|\Delta(g_{k_1})|$,
\begin{align*}
\omega_2(h_1,\dots,h_{k_2})&=
\omega_2(z^>_j(g_{k_1}) h_1,\dots,h_{k_2})\\&-
\sum_{i=1}^{k_2-1}(-1)^{i}\omega_2(z^>_j(g_{k_1}),h_1,\dots,h_i h_{i+1},\dots,h_{k_2})\\&-
(-1)^{k_2}\omega_2(z^>_j(g_{k_1}),h_1,\dots,h_{k_2-1}).
\end{align*}
We use the first relation to expand $(\omega_1\cup\eta_2)(g_1,\dots,h_{k_2})$:
\begin{align*}
&(\omega_1\cup\eta_2)(g_1,\dots,h_{k_2})\\=&
\omega_1(g_1,\dots,g_{k_1})\cdot\eta_2(h_1,\dots,h_{k_2})\\=&
\omega_1(g_1,\dots,g_{k_1})\cdot
\sum_{j=1}^{|\Delta(h_1)|}
\varphi(h_1^{(j)})\cdot\omega_2(z^>_j(h_1),h_2,\dots,h_{k_2})\\=&
\sum_{j=1}^{|\Delta(h_1)|}
\omega_1(g_1,\dots,g_{k_1}z^<_j(h_1))\cdot
\varphi(h_1^{(j)})\cdot\omega_2(z^>_j(h_1),h_2,\dots,h_{k_2})\\+&
\sum_{j=1}^{|\Delta(h_1)|}
(-1)^{k_1}\sum_{i=1}^{k_1-1}(-1)^{i}\omega_1(g_1,\dots,g_i g_{i+1},\dots,g_{k_1},z^<_j(h_1))\cdot
\varphi(h_1^{(j)})\cdot\omega_2(z^>_j(h_1),h_2,\dots,h_{k_2})\\-&
\sum_{j=1}^{|\Delta(h_1)|}
(-1)^{k_1}\omega_1(g_2,\dots,g_{k_1},z^<_j(h_1))\cdot
\varphi(h_1^{(j)})\cdot\omega_2(z^>_j(h_1),h_2,\dots,h_{k_2}).
\end{align*}
In the same way, we use the second relation to expand $(\eta_1\cup\omega_2)(g_1,\dots,h_{k_2})$:
\begin{align*}
&(\eta_1\cup\omega_2)(g_1,\dots,h_{k_2})\\=&
\eta_1(g_1,\dots,g_{k_1})\cdot\omega_2(h_1,\dots,h_{k_2})\\=&
\sum_{j=1}^{|\Delta(g_{k_1})|}
\omega_1(g_1,\dots,g_{k_1-1},z^<_j(g_{k_1}))\cdot\varphi(g_{k_1}^{(j)})
\cdot\omega_2(h_1,\dots,h_{k_2})\\=&
\sum_{j=1}^{|\Delta(g_{k_1})|}
\omega_1(g_1,\dots,g_{k_1-1},z^<_j(g_{k_1}))\cdot\varphi(g_{k_1}^{(j)})\cdot
\omega_2(z^>_j(g_{k_1}) h_1,\dots,h_{k_2})\\-&
\sum_{j=1}^{|\Delta(g_{k_1})|}
\omega_1(g_1,\dots,g_{k_1-1},z^<_j(g_{k_1}))\cdot\varphi(g_{k_1}^{(j)})\cdot
\sum_{i=1}^{k_2-1}(-1)^{i}\omega_2(z^>_j(g_{k_1}),h_1,\dots,h_i h_{i+1},\dots,h_{k_2})\\-&
\sum_{j=1}^{|\Delta(g_{k_1})|}
\omega_1(g_1,\dots,g_{k_1-1},z^<_j(g_{k_1}))\cdot\varphi(g_{k_1}^{(j)})\cdot
(-1)^{k_2}\omega_2(z^>_j(g_{k_1}),h_1,\dots,h_{k_2-1}).
\end{align*}
Finally, we define
$$
\eta(g_1,\dots,g_{k_1-1},e,h_2,\dots,h_{k_2})=
\sum_{j=1}^{|\Delta(e)|}\omega_1(g_1,\dots,g_{k_1-1},z^<_j(e))
\cdot\varphi(e^{(j)})\cdot
\omega_2(z^>_j(e),h_2,\dots,h_{k_2})
$$
and write out $\delta\eta(g_1,\dots,g_{k_1},h_1,\dots,h_{k_2})$:
\begin{align*}
&
\sum_{j=1}^{|\Delta(h_1)|}\omega_1(g_2,\dots,g_{k_1},z^<_j(h_1))
\cdot\varphi(h_1^{(j)})\cdot
\omega_2(z^>_j(h_1),h_2,\dots,h_{k_2})\\+&
\sum_{i=1}^{k_1-1}(-1)^i\sum_{j=1}^{|\Delta(h_1)|}
\omega_1(g_1,\dots,g_i g_{i+1},\dots,g_{k_1},z^<_j(h_1))
\cdot\varphi(h_1^{(j)})\cdot
\omega_2(z^>_j(h_1),h_2,\dots,h_{k_2})\\+&
(-1)^{k_1}\sum_{j=1}^{|\Delta(g_{k_1}h_1)|}
\omega_1(g_1,\dots,z^<_j(g_{k_1} h_1))
\cdot\varphi((g_{k_1}h_1)^{(j)})\cdot
\omega_2(z^>_j(g_{k_1} h_1),h_2,\dots,h_{k_2})\\+&
\sum_{i=1}^{k_2-1}(-1)^{k_1+i}\sum_{j=1}^{|\Delta(g_{k_1})|}
\omega_1(g_1,\dots,g_{k_1-1},z^<_j(g_{k_1}))
\cdot\varphi(g_{k_1}^{(j)})\cdot
\omega_2(z^>_j(g_{k_1}),h_1,\dots,h_i h_{i+1},\dots,h_{k_2})\\+&(-1)^{k_1+k_2}
\sum_{j=1}^{|\Delta(g_{k_1})|}
\omega_1(g_1,\dots,g_{k_1-1},z^<_j(g_{k_1}))
\cdot\varphi(g_{k_1}^{(j)})\cdot
\omega_2(z^>_j(g_{k_1}),h_1,\dots,h_{k_2-1}).
\end{align*}
Now most of the summands of $(\omega_1\cup\eta_2+\eta_1\cup\omega_2-(-1)^{k_1}\delta\eta)(g_1,\dots,g_{k_1},h_1,\dots,h_{k_2})$ cancel out and we are left with
\begin{align*}
&
\sum_{j=1}^{|\Delta(g_{k_1})|}
\omega_1(g_1,\dots,g_{k_1-1},z^<_j(g_{k_1}))\cdot\varphi(g_{k_1}^{(j)})\cdot
\omega_2(z^>_j(g_{k_1}) h_1,\dots,h_{k_2})\\+&
\sum_{j=1}^{|\Delta(h_1)|}
\omega_1(g_1,\dots,g_{k_1}z^<_j(h_1))\cdot
\varphi(h_1^{(j)})\cdot\omega_2(z^>_j(h_1),h_2,\dots,h_{k_2})\\-&
\sum_{j=1}^{|\Delta(g_{k_1}h_1)|}
\omega_1(g_1,\dots,z^<_j(g_{k_1} h_1))
\cdot\varphi((g_{k_1}h_1)^{(j)})\cdot
\omega_2(z^>_j(g_{k_1} h_1),h_2,\dots,h_{k_2}).
\end{align*}
Each of the three sums above corresponds to a side of the geodesic triangle $(1,g_{k_1},g_{k_1}h_1)$. Using property (iv) of decompositions for $g_{k_1}$ and $h_1$ we are given the elements $c_1$, $c_2$, $c_3$, $r_1$, $r_2$, $r_3\in F$.
Since $(g_1,\dots,g_{k_1},h_1,\dots,h_{k_2})$ is in $\mathcal{B}^{k_1+k_2}$, $g_{k_1}h_1$ is a reduced word and $c_2$ is the neutral element. Furthermore, for $j\leq|\Delta(c_1)|$, we have $z^<_j(g_{k_1})=z^<_j(g_{k_1} h_1)$ and $z^>_j(g_{k_1}) h_1=z^>_j(g_{k_1} h_1)$, hence the first $|\Delta(c_1)|$ terms of the first sum cancel with the first $|\Delta(c_1)|$ terms of the third sum. A similar observation can be done with the last $|\Delta(c_3)|$ terms of the second sum and the last $|\Delta(c_3)|$ terms of the third sum. Therefore, only the summands corresponding to $\Delta(r_1)$, $\Delta(r_2)$, $\Delta(r_3)$ survive, which means that there are at most $3R$ terms left. Each of these terms is bounded because $\omega_1$ and $\omega_2$ are bounded and $\varphi$ is bounded on the pieces of the $\Delta$-decomposition.
\begin{figure}[H]
    \centering
    \includegraphics[scale=1.2]{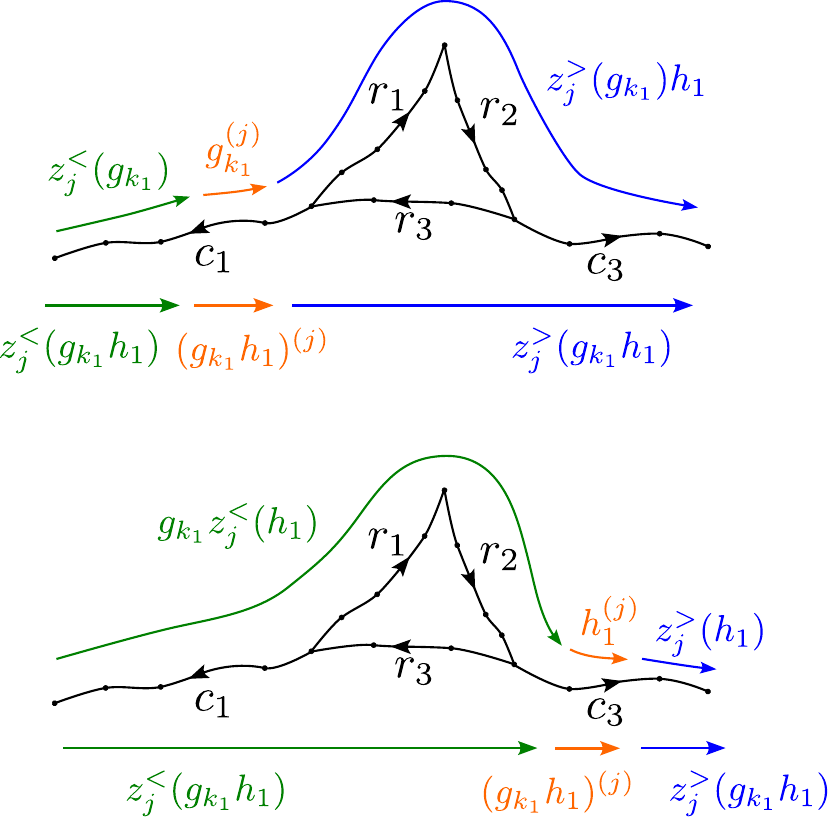}
    \caption{Cancellation of the first $|\Delta(c_1)|$ terms (above) and of the last $|\Delta(c_3)|$ terms (below).}
\end{figure}
\section{Proof of Theorem 2}
Let $M$ be a negatively curved orientable compact Riemannian manifold with (possibly empty) convex boundary. The universal covering $\widetilde{M}$ is continuously uniquely geodesic and thus for every $(x_0,\dots,x_k)\in\widetilde{M}^{k+1}$, by repeatedly coning on the $x_i$ one can define the straight $k$-simplex $[x_0,\dots,x_k]\in\S_k(\widetilde{M})$ as constructed in Section 8.4 of \cite{frigerio2017bounded}. 
The fundamental group $\pi_1(M)=\Gamma$ acts on the universal covering $\widetilde{M}$ via deck transformations and this defines in turn an action of $\Gamma$ on $C^{\bullet}_b(\widetilde{M})$. We denote by $\C_b^\bullet(\widetilde{M})^\Gamma$ the subcomplex of $\Gamma$-invariant cochains.
The covering map $p\colon\widetilde{M}\to M$ induces an isometric isomorphism of normed complexes $\C_b^{\bullet}(M)\xrightarrow[]{\cong}\C_b^{\bullet}(\widetilde{M})^{\Gamma}$.
Similarly, $\Gamma$ acts on $\Omega^k(\widetilde{M})$, and we denote by $\Omega^k(\widetilde{M})^\Gamma$ the space of $\Gamma$-invariant $k$-forms of $\widetilde{M}$. By pulling-back via the covering projection we get the identification $\Omega^k(M)\xrightarrow[]{\cong}\Omega^k(\widetilde{M})^\Gamma$.
\begin{notation}
We will define cochains only on straight simplices with the understanding that they can be extended to every simplex $s\in\S_k(\widetilde{M})$ by precomposing with the \emph{straightening operator} $s\mapsto [s(e_0),\dots,s(e_k)]$, where $e_0,\dots,e_k$ are the vertices of the standard simplex $\Delta^k$ (see e.g. Section 8.7 of \cite{frigerio2017bounded}).
\end{notation}
For any $\psi\in\Omega^k(\widetilde{M})^\Gamma$, we define a cochain $c_\psi\in\C^k(\widetilde{M})^{\Gamma}$ by setting for every $(x_0,\dots,x_k)\in\widetilde{M}^{k+1}$,
$$
c_\psi(s)=\int_{[x_0,\dots,x_k]}\psi.
$$
Let $\alpha_1\in\H_b^{k_1}(M)$ and $\alpha_2\in\H_b^{k_2}(M)$. It is shown in Section 2.3 of \cite{https://doi.org/10.48550/arxiv.2202.04419} that there are representatives $\omega_i\in\C_b^{k_i}(\widetilde{M})^\Gamma$ for $\alpha_i$ that depend smoothly on the vertices of simplices. More precisely, there are bounded smooth functions $f_{\omega_i}\colon\widetilde{M}^{k_i+1}\to\mathbb{R}$ such that for every $(x_0,\dots,x_k)\in\widetilde{M}^{k+1}$, $$\omega_i([x_0,\dots,x_{k_i}])=f_{\omega_i}(x_0,\dots,x_k).$$
Let $\xi=d\varphi\in\Omega^2(M)$ be an exact 2-form. It is proved in \cite{https://doi.org/10.48550/arxiv.2202.04419} that $[c_\xi]\cup\alpha_i=0\in\H_b^{k_i+2}(M)$ and thus $\langle\alpha_1,[c_\xi],\alpha_2\rangle$ is well defined.
The bounded primitives $\beta_1\in\C_b^{k_1+1}(\widetilde{M})^\Gamma$ and $\beta_2\in\C_b^{k_2+1}(\widetilde{M})^\Gamma$ of the cup products
$\omega_1\cup \delta\varphi=\delta\beta_1$ and $\delta\varphi\cup\omega_2=\delta\beta_2$ are given by
$$
\beta_1=(-1)^{k_1}\omega_1\cup c_\varphi-\delta\eta_1
$$
$$
\beta_2=c_\varphi\cup\omega_2+\delta\eta_2,
$$
where $\eta_1\in\C^{k_1}(F)$ and $\eta_2\in\C^{k_2}(F)$ are defined as follows:
\begin{align*}
\eta_1([x_0,\dots,x_{k_1}])&=\int_{[x_{k_1-1},x_{k_1}]}f_{\omega_1}(x_0,\dots,x_{k_1-1},-)\cdot\varphi\\
\eta_2([x_0,\dots,x_{k_2}])&=\int_{[x_0,x_1]}\varphi\cdot f_{\omega_2}(-,x_1,\dots,x_{k_2}).
\end{align*}
Exactly like in the case of the free groups we can now use these primitives to write out the following representative of an element in $\langle\alpha_1,[c_{d\varphi}],\alpha_2\rangle$:
\begin{align*}
(-1)^{k_1}\omega_1\cup\beta_2-\beta_1\cup\omega_2&=
(-1)^{k_1}\omega_1\cup(c_\varphi\cup\omega_2+\delta\eta_2)
-((-1)^{k_1}\omega_1\cup c_\varphi-\delta\eta_1)\cup\omega_2
\\&=(-1)^{k_1}\omega_1\cup\delta\eta_2+\delta\eta_1\cup\omega_2.
\end{align*}
The cochain $\omega_1\cup\eta_2+\eta_1\cup\omega_2\in\C^{k_1+k_2}(\widetilde{M})^\Gamma$ is a primitive for the cocycle above, but it is not necessarily bounded. In order to find a bounded primitive we proceed as we did in the last section. We define as follows an element $\eta\in\C^{k_1+k_2-1}(\widetilde{M})^\Gamma$ such that $\omega_1\cup\eta_2+\eta_1\cup\omega_2-(-1)^{k_1}\delta\eta$ is bounded:
$$
\eta([x_0,\dots,x_{k_1+k_2-1}])=\int_{[x_{k_1-1},x_{k_1}]}f_{\omega_1}(x_0,\dots,x_{k_1-1},-)\cdot\varphi\cdot f_{\omega_2}(-,x_{k_1},\dots,x_{k_1+k_2-1}).
$$
In fact, after some computations (we will skip some details here since they are completely analogous -and easier- to the computations done in the previous section) we get that
\begin{align*}
&(f_{\omega_1}\cup\eta_2+\eta_1\cup f_{\omega_2}-(-1)^{k_1}\delta\eta)([x_0,\dots,x_{k_1+k_2}])
\\=&
\int_{[x_{k_1-1},x_{k_1}]}
f_{\omega_1}(x_0,\dots,x_{k_1-1},-)
\cdot\varphi\cdot
f_{\omega_2}(-,x_{k_1+1},\dots,x_{k_1+k_2})\\+&
\int_{[x_{k_1},x_{k_1+1}]}
f_{\omega_1}(x_0,\dots,x_{k_1-1},-)
\cdot\varphi\cdot f_{\omega_2}(-,x_{k_1+1},\dots,x_{k_1+k_2})\\-&
\int_{[x_{k_1-1},x_{k_1+1}]}
f_{\omega_1}(x_0,\dots,x_{k_1-1},-)
\cdot\varphi\cdot
f_{\omega_2}(-,x_{k_1+1},\dots,x_{k_1+k_2})
\\=&
\int_{\partial[x_{k_1-1},x_{k_1},x_{k_1+1}]}
f_{\omega_1}(x_0,\dots,x_{k_1-1},-)
\cdot\varphi\cdot
f_{\omega_2}(-,x_{k_1+1},\dots,x_{k_1+k_2})
\\=&
\int_{[x_{k_1-1},x_{k_1},x_{k_1+1}]}
d(f_{\omega_1}(x_0,\dots,x_{k_1-1},-)
\cdot\varphi\cdot
f_{\omega_2}(-,x_{k_1+1},\dots,x_{k_1+k_2})).
\end{align*}
The integration domain is now $[x_{k_1-1},x_{k_1},x_{k_1+1}]$, which is a straight triangle with bounded area. Therefore, to conclude the proof it is enough to check that the differential form $d(f_{\omega_1}(x_0,\dots,x_{k_1-1},-)
\cdot\varphi\cdot
f_{\omega_2}(-,x_{k_1+1},\dots,x_{k_1+k_2}))$ assumes uniformly bounded values over orthonormal frames. This is straightforward to see by using that $M$ is compact and using the bounds on $df_{\omega_i}$ stated in Lemma 1 of \cite{https://doi.org/10.48550/arxiv.2202.04419}.
\printbibliography
\end{document}